\documentclass[referee,sn-basic]{sn-jnl}

\usepackage{graphicx}%
\usepackage{multirow}%
\usepackage{amsmath,amssymb,amsfonts}%
\usepackage{amsthm}%
\usepackage{mathrsfs}%
\usepackage[title]{appendix}%
\usepackage{xcolor}%
\usepackage{textcomp}%
\usepackage{manyfoot}%
\usepackage{booktabs}%
\usepackage{algorithm}%
\usepackage{algorithmicx}%
\usepackage{algpseudocode}%
\usepackage{listings}
\usepackage{stackengine}

\theoremstyle{thmstyleone}%
\newtheorem{proposition}{Proposition}
\newtheorem{lemma}
{Lemma}
\theoremstyle{thmstyletwo}%
\newtheorem{remark}{Remark}%

\theoremstyle{thmstylethree}%
\newtheorem{definition}{Definition}%

\begin{document}

\title[Article Title]{A Note on Inferential Decisions, Errors and Path-Dependence}


\author*[1]{\fnm{Ken Kangda} \sur{Wren}}\email{k.k.wren@york.ac.uk}



\affil*[1]{\orgdiv{Department of Mathematics}, \orgname{University of York}, \orgaddress{\street{Heslington}, \city{York}, \postcode{YO10 5DD}, 
\country{UK}}}




\abstract{Consider sequential binary testing under model uncertainty or misspecification in an otherwise 'ideal' setting: the \emph{a posteriori} belief process and its objective conditional probability counterpart may differ but converge to the same correct outcome. We show that under common conditions (defined) unless the two are 'essentially identical', differing only by \emph{a priori} factors, time-homogeneous continuous decisions based on one must fail to be path-independent with respect to state-variables based on the other or any non-essentially-identical process. The difference between them, inferential error, decomposes uniquely into two independent components: a path-independent systematic bias and a path-dependent, not necessarily systematic, error.}




\maketitle
\newpage\section*{Declarations}

\begin{itemize}
\item \textbf{Author K.K.Wren} declares:
\item This research has received no funding.
\item I have no competing interests to disclose.
\item I approve the Ethics, Consent to Participate, and Consent to Publish Declarations.
\item There is no Data availability issue to address in relation to this manuscript
\item I wrote and reviewed the manuscript. No other authors are involved.
\end{itemize}\newpage
\section{Introduction}\label{sec1}
Decision under uncertainty and noisy data is formulated routinely with inferential belief (or equivalent) at the core, to be fed into some algorithm, for 'optimal' outputs: thus raw data are distilled into 'rational beliefs' and the select algorithm takes care of 'statistical noise'. 
Concerns have been raised in recent years about this framework, in particular, the effects of model uncertainty or misspecification and how to address them (e.g. \citet{BB08}, \citet{JCB}, \citet{lyl}). These suggest that the sequence of modelling and testing choices by which evidence is generated may have to be taken into account; see \citet{BK21} for further reading on the need for a broader framework for statistical evaluation.

We report a finding in line with the above: for binary inference under misspecified models and common conditions, inferential error, of a given belief process vs another (e.g. true probability), is in general not {path-independent} (not round-tripping despite the given belief process doing so); it decomposes uniquely into two \emph{independent} components, one {path-independent}, systematic  and identifiable as a well-known behavioural bias, one path-dependent and not necessarily systematic.

A novel aspect of the report is the isolation and characterisation of the effects of model misspecification in terms of 
how systems driven by statistical inference may be modelled. Such systems, mostly {socio-economic} and operating {without} 'the best and the latest' statistical know-how (e.g. institution, economy, market), contend with real-world model uncertainty and misspecification, the effects of which must be characterised for the better modelling and understanding of these important systems
. For instance, our result affects asset-market modelling\footnote{A 'market' is a system driven by an objective conditional probability process, one reference inference process, either the same as the first ('orthodox theory') or different ('reality'), and another inference process '{equivalent}' to the reference (to attain 'viable pricing'); see {\citet{coh}} for background
.} in two interesting ways, both via path-in/dependence. %
These are analogous to, but more structural and tangible than, \citet{GTim}'s observation that learning (a missing continuous parameter) makes asset-prices depend on the sample-path of data.

For this study, to be clear, we take the usual setup of sequential binary testing \emph{except} that the test models can be misspecified. In this controlled setting and for systems with variables adapted to the tests, path-in/dependence 
and 'inferential redundancy' are defined and examined, culminating in the main results, examples and implications, including the {decomposition} of inferential error into a 
path-independent systematic bias for/against 'the status quo' and a path-dependent error that may or may not be systematic (although it is in typical examples under i.i.d data).

\section{Setup}
\subsection{Path-Independence in a Simple Backdrop}
Consider {data process} \(\{{D}_n\}\), with natural filtration \(\{\mathcal{F}_n\}\) and 
law \(\mathbf{P}_b\) whose identity is in question, on filtered space \((\Omega,\{\mathcal{F}_n\}; \mathbf{P}_b)\), where \emph{cumulative} data \(D_n:\Omega
\equiv\mathbf{V}^{\mathbb{N}}\ni\omega
\mapsto D_n(\omega)\ni\mathbf{V}^n
\) is a \(n\)-string of real vectors \(D_n[i](\omega)\in\mathbf{V}\), \(i=1,...,n\in\mathbb{N}\).

Raw data being usually 'too noisy', it is desirable to model a system driven by it through an array \(\{X_n\}\) of \(\{\mathcal{F}_n\}\)-adapted
state-variables that capture the system-relevant features of the data. By this goal, any dependence beyond 'the current state' signals some deficiency of the model devised or chosen.

We shall focus on systems having: 1) a state-variable that is an 
inference process \(\{\pi^{\theta}_n\}\) about an outcome \(B\in\{\theta,\theta_o\}\), where \(\theta_o\) is some \emph{baseline} (e.g. 'status quo', 'no effect') and 2) key properties that are \emph{time-homogeneous continuous functions} of this or another inference process about \(B\) based on the same data.
\begin{definition}\label{PIdefo}
         A process \(\{Y_n\}\) is \emph{path-independent} with respect to \(\{\mathcal{F}_n\}\)-adapted process \(\{X_n\}\) satisfying \(Supp(X_n)\supseteq Supp(X_{n'})\), \(\forall n>n'\in\mathbb{N}\) finite, if: 1) \(X_n(\omega)=X_{n'}(\omega')\implies Y_n(\omega)=Y_{n'}(\omega')\), \(\forall \omega,\omega'\in\Omega\) almost surely; 2) \(\sigma({Y_n})\subseteq\sigma({X_n})\), \(\forall n\in\mathbb{N}\).
\end{definition}
\begin{remark}\label{strict}
'Path-independence' often means point-1 ('round-tripping'), saying little about the underlying relationship. Point-2 makes \(Y_n\) a function \(f_n\) 
    of \(X_n\) 
    (Proposition 4.9, \citet{4.9}); its fail says 'dependence beyond \(X_n\)'; its pass implies point-1 iff. \(f_n\equiv f_{n'}\). \end{remark}
    
    Violation of our strict notion means dependence on \(\mathcal{F}^X_n\setminus\sigma(X_n)\) (state history), on \(\mathcal{F}_n\setminus\sigma(X_n)\) (data history), or on 
time (passage of history)
. The last one is arguably not 'real path-dependence'. This is moot, as we limit ourselves to {time-homogeneous dependence} on inference throughout; see also Remark \ref{clarify'} later. 
    \subsection{Common Sequential Inference}\label{st}
Consider the sequential testing of \(B\in\{\theta,\theta_o\}\) using data \(\{{D}_n\}\) and \emph{test measure} \(Q_\theta\) (test model) vs \emph{baseline measure} \(Q_{{\theta_o}}\) (null model). Write the objective law of data given \(B=b\) as \(\mathbf{P}_{b}\)
. 
The following \emph{common test features} apply throughout:
    \begin{enumerate}
    \item\label{notime}All the measures (models/laws) are time-homogeneous\footnote{There thus can be no time-dependence in the inferential dynamic other than that from data accumulation or any time-varying parameters of the dataflow itself.};
    \item\label{uni}The baseline measure is objectively accurate ('universal'): \(Q_{{\theta_o}}\equiv \mathbf{P}_{{\theta_o}}\);    
    \item{\label{simsim} On partial-data space \((\mathbf{V}^{n}, \mathcal{F}_n)\), \(n\in\mathbb{N}\), {equivalence} \(Q_\theta|_{n}\sim Q_{{\theta_o}}|_{n}\sim \mathbf{P}_{\theta}|_{n}\) holds;}
    \item{\label{regtest} On total-data space \((\mathbf{V}^{\mathbb{N}}, \mathcal{F}_\infty)\) mutual singularity \(Q_\theta\perp Q_{{\theta_o}}\) and \(\mathbf{P}_{\theta}\perp \mathbf{P}_{\theta_o}\) hold, with 
    \(\mathbf{P}_{\theta}\prec\prec Q_{\theta}\), ensuring almost sure, objectively accurate, \(B\)-detection at \(\infty\).}\end{enumerate}
    
    Testing is based on \emph{likelihood ratio} (LR) \(L^\theta_n(\cdot):=\frac{dQ_\theta|_n(\cdot)}{dQ_{{\theta_o}}|_{n}(\cdot)}\), \(n\in\mathbb{N}\), whose conditional versions \(L^{\theta}_{n|m}(\cdot):=\frac{dQ_\theta|_{n}(\cdot|\mathcal{F}_{m})}{dQ_{{\theta_o}}|_{n}(\cdot|\mathcal{F}_{m})}\), \(m<n\), are assumed to exist
. Then: 
\begin{align}
    \label{the1}&{E}_{{o}}[L^\theta_n]\equiv1\equiv{E}_{{o}}[L^\theta_{n|m}|\mathcal{F}_{m}],\text{ }\forall n>m\in\mathbb{N},\\
    \label{lr}&L^{\theta}_n(\cdot|\mathcal{F}_{m})\equiv
L^{\theta}_{n|m}(\cdot)L^{\theta}_m,\text{ }\forall n>m\in\mathbb{N},\\
\label{resolving}&l^{\theta}_n(B):=\log L^{\theta}_n(B)\to(-1)^{\mathbf{1}_{\{B={\theta_o}\}}}\cdot\infty\text{ only as \({n\to\infty}\)};    
\end{align}
where the \emph{log-LR process} \(\{l^{\theta}_n\}\) best describes the test dynamic (see Appendix \ref{infps}). 

Any LR level \(L^{\theta}_n\in\mathbb{R}^+\) becomes some \emph{a posteriori} belief \({\pi}^\theta_n\) given any \emph{a priori} belief \({\pi}^{\theta}_0\in(0,1)\), via Bayes' Rule, written in \emph{inferential odds} \(L[\pi]:=\frac{\pi}{1-\pi}
\)
,\begin{equation}\label{br1}
    L[{\pi}^{\theta}_n(\cdot|\mathcal{F}_{m})]= L[{\pi}^{\theta}_0]\times 
    {L^{\theta}_n(\cdot|\mathcal{F}_{m})}
    =L[{\pi}^{\theta}_m]\times 
    {L^{\theta}_{n|m}(\cdot)},\text{ }\forall n>m\in\mathbb{N}.
    \end{equation}Note the independence of \(L[{\pi}^{\theta}_0]\) and \(\{L^{\theta}_n\}\), and of \(L[{\pi}^{\theta}_0]\) and \(\{L^{\theta}_{n|m}\}\), by definition.
    
    \emph{Common inferential \(B\)-testing} occurs on \(
    \{\theta,{\theta_o}\}\times\Omega\ni(B,\omega)\), under the natural measure \(\mathbb{Q}:={\pi}^{B
    }_0\times Q_{B
    }\), \({\pi}^{B}_0\in(0,1)\); note \(\pi^{{\theta}}_n\equiv\mathbb{Q}(\{B={\theta}\}|\mathcal{F}_n)\), \(\forall n\in\mathbb{N}\).

\subsection{Inferential Redundancy}
Consider two common tests, based on data \(\{D_n\}\), one using test model \(Q_\theta
\) and the other, \(\hat{Q}_\theta
\), generating respective LR process \(\{L^{\theta}_n\}\) and \(\{\hat{L}^{\theta}_n\}\), associated with which are the \emph{a posteriori} beliefs \(\{\pi^\theta_n\}\) and \(\{\hat{\pi}^\theta_n\}\), given \emph{a priori} belief \({\pi}^\theta_{0}\) and \(\hat{\pi}^\theta_0\), respectively.

\begin{definition}\label{d1}
    A common test, with test model \(\hat{Q}_\theta\) and 
    \emph{a posteriori} beliefs \(\{\hat{\pi}^\theta_n\}\), is said to be \emph{inferentially redundant} with respect to another, with 
     \({Q}_\theta\) and \(\{{\pi}^\theta_n\}\), if: 
     \(\sigma({\hat{\pi}^\theta_n})\subseteq\sigma({\pi^\theta_n})\), \(\forall n\in\mathbb{N}\), such that the implied dependence is continuous\footnote{
    Continuity makes 'small evidence' remain 'small'. Note also that although the definition is asymmetric, redundancy, as will be seen, is always mutual.}.
     
\end{definition}

\subsection{Sufficiently Rich Tests}\label{technical}
Both path-independence and inferential redundancy concern links between paths, not likelihoods. Define thus \emph{rich tests}, where finite conditioning never permanently lowers 'future abundance' (e.g. it may affect only the shapes, not the supports, of distributions, as in the 'typical tests' of Remark \ref{kaku} %
 below).
 
 To be exact, let the following hold at any \(m\in\mathbb{N}\) finite, irrespective of conditioning by outcomes of the form \(\{L^{\theta}_m=(\cdot)\in\mathbb{R}^+\}\):
    \begin{enumerate}
        \item\label{p.1} 
    conditional-LR variables \(L^{\theta}_{n|m}\) given \(L^{\theta}_m=(\cdot)\), \(m<n\), have \emph{conditional supports} \(\mathcal{L}^{(\cdot)}_{n|m}\) 
    that satisfy the following: for any given pair of \(m\&m'\in\mathbb{N}\) and conditioning outcomes \((\cdot)\) and \((\cdot)'\), \(\exists N\in\mathbb{N}\) finite 
        such that,
        \begin{align}
       \label{spam'''}&\mathcal{L}^{(\cdot)}_{n|m}\cap\mathcal{L}^{(\cdot)'}_{n|m'}\not=\emptyset,\text{ \(\forall n>N\)};
\end{align}
        \item \label{p2}\emph{total conditional-support} \(\mathcal{L}^{(\cdot)}_{|m}:=\bigcup_{n>m}^{\infty}
        \mathcal{L}^{(\cdot)}_{n|m}
        \) 
        satisfies: \(\forall m\in\mathbb{N}\) finite,
        \begin{equation}\label{denseness}
        {\mathcal{L}^{(\cdot)}_{|m}}\subset_{dense}
        (\inf\mathcal{L}^{(\cdot)}_{|m}, \sup\mathcal{L}^{(\cdot)}_{|m})=\mathbb{R}^+
        .
        \end{equation}
   
    \end{enumerate}
    
    The term '\emph{common}' (tests) shall signify the default pre-condition, namely, all the properties given in Section \ref{st} and 
\ref{technical},
for our arguments and conclusions.
\begin{remark}\label{kaku}
    'Typical tests' draw i.i.d data from a distribution on \(\mathbb{R}^k\), \(\mu_{\theta_o}\) or \(\mu_{{\theta}}\sim\mu_{\theta_o}\); the limit measures \(\mu^{\mathbb{N}}_{\theta}\) and \(\mu^{\mathbb{N}}_{{\theta_o}}\) on \((\mathbb{R}^k)^{\mathbb{N}}\) exist, with \(\mu^{\mathbb{N}}_{\theta}\perp\mu^{\mathbb{N}}_{{\theta_o}}\) unless \(\mu_{\theta}=\mu_{{\theta_o}}\) (Kakutani's Theorem). When \(\mu_{\theta}\) and \(\mu_{{\theta_o}}\) have full support, such tests are trivially 'rich'. Even when they are binary (e.g. coinflips), property (\ref{spam'''}-\ref{denseness}) are easily verified. 
\end{remark}

\section{Main Results and Proof}\label{lemmasec}

\begin{lemma}\label{l1}
    Inferential redundancy between common tests about the same outcome using the same data must take the form below: 
    for some \(\gamma\in(0,1]\) and \(\{c_n\}>0\), all data-independent,
    \begin{align}\label{LRitself'}
        &\hat{L}^{\theta}_{n}=(L^{\theta}_{n})^{\gamma}c_n,
        \text{ }\hat{L}^{\theta}_{n|m}=(L^{\theta}_{n|m})^{\gamma}\frac{c_n}{c_m},\text{ \(
        \forall m<n\in\mathbb{N}\), subject to:}\\\label{normlstn}&{E}_{{o}}[(L^{\theta}_{n})^\gamma]\equiv\frac{1}{c_n}\le1\text{ }\&\text{ }
    {E}_{{o}}[(L^{\theta}_{n|m})^\gamma|\mathcal{F}_m]\equiv\frac{c_m}{c_n}\le1.
\end{align}
\end{lemma}

\begin{proof}
Redundancy maps obey inference equations that become Cauchy under 'richness' and continuity; 
  the solutions are positive power-laws (Appendix \ref{proofl1}). Labelling choice makes \(\gamma\le1\); normalisation (\ref{normlstn}) stems from (\ref{the1}), and the bounds of \(\{c_n\}\), Jensen's inequality.
\end{proof}

    Note that \(\gamma=1\) implies \(c_n\equiv1\) ((\ref{normlstn})); so \(\hat{L}^\theta_n\equiv L^\theta_n\). 
    Conversely, path-independent inferential redundancy demands \(c_n\equiv c_{n'}\) (Definition \ref{PIdefo}); so \(\gamma=1\), by (\ref{LRitself'}-\ref{normlstn}), Jensen's inequality and the strict concavity of power-law with \(\gamma\in(0,1)\). Hence:
\begin{proposition}\label{p1}
    No common inferential test can be path-independent to another about the same outcome using the same data without being identical to it up to {a priori} factors.
\end{proposition}

\begin{remark}\label{clarify}
      All above apply to tests without shared null-models (Remark \ref{prfrmrk}, Appendix \ref{proofl1}).
\end{remark}
\begin{remark}\label{clarify'}
Redundancy with time-dependence (\(\gamma\not=1\)) exists, for 'just so' pairings. In reality \(\hat{Q}_\theta\not=Q_\theta\) all but incurs genuine (state/data) path-dependence. Consider the exception: two 'typical tests' (Remark \ref{kaku}) of Gaussians of differing means: \(\mu_{\theta}\)-vs-\(0\) \& \(\hat{\mu}_{{\theta}}\)-vs-\(0
     \). Let the second be objectively accurate, driving the conditional probability of outcome \(\theta=\{cancer!\}\). Say \(\hat{\mu}_{{\theta}}<{\mu}_{{\theta}}\), then \(\hat{L}^\theta_{n|m}=(L^\theta_{n|m})^{\hat{\mu}_{\theta}/{\mu_{\theta}}}e^{(n-m)\cdot\hat{\mu}_{{\theta}}({\mu_{\theta}-\hat{\mu}_{{\theta}}})}\) (cf. (\ref{LRitself'})). 
     This, with true mean \(\hat{\mu}_{{\theta}}\) unknown, 
     is no less challenging in practice than genuine path-dependence, especially when 'absolute time' (\(m\) or \(n\) itself, not just \(n-m\)) is often ambiguous outside of controlled experiments. If the tests are coinflips, where \(\mu_{\theta}\) (\(\hat{\mu}_{{\theta}}\)) are heads-probabilities, and '\(0\)', the fair-coin baseline, no such or any other redundancy is possible: all coinflip-tests with \(\mu_{\theta}\not=\hat{\mu}_{{\theta}}\) are non-redundant and mutually path-dependent.
\end{remark}

\section{Implications for Inference-adapted Dynamics}\label{PvQ}
Often the objective conditional probability process \(\{p
^\theta_n\}\) is known up to equivalence only: some \(\{\hat{\pi}^\theta_n\}\sim\{p
    ^\theta_n\}
    \) given identical data. By Proposition \ref{p1}, 
unless they differ only by \emph{a priori} factors, system properties of the form \(\{Y^p_n\}:=\{u[{p^\theta_n}]\}\) cannot be path-independent to state-variables of the form \(\{{X}^{\hat{\pi}}_n\}:=\{v[\hat{\pi}^\theta_n]\}\), where \(u\) and \(v\) are continuous functions. This is true of any inference \(\{{\pi}^{\theta}_n\}\sim\{\hat{\pi}^{\theta}_n\}\) under identical data. For instance, properties driven by 'difference of opinions', \(\{X^{\hat{\pi}}_n\}-\{X^{\pi}_n\}\), must be path-dependent to both 
unless \(\{{\pi}^{\theta}_n\}\) and \(\{\hat{\pi}^{\theta}_n\}\) differ only by \emph{a priori} factors. 

\begin{remark}\label{asspximp}
    This has a quick implication for asset-market modelling, by which an asset under binary risk \(B\in\{\theta, \theta_o\}\) is priced by some inference process \(\{\pi^\theta_n\}
    \) that is a so-called Risk-Neutral Equivalent (RNE) law of \(B\) given true law \(\{p^\theta_n
    \}
    \); the reward for bearing \(B\)-risk, \(expected\text{ }income\)\(-\)\(price\text{ }paid\), has the form \(\{Y^{p
    }_n\}-\{Y^\pi_n\}\); it is \emph{viable} ('arbitrage-free') iff. \(\{\pi^\theta_n\}\sim\{p
    ^\theta_n\}
    \)
. However, economic theories dictate that price \(Y^\pi_n\) maximises expected agent-utility, thereby making \(\pi^\theta_n\) a time-homogeneous function of \(p
    ^{\theta}_n\). Proposition \ref{p1} then implies 
    that \emph{economic asset-pricing} has \emph{one} viable form (\emph{one} type of RNE law) only. Further, as true law \(\{p^\theta_n
    \}\) is often unknown \emph{ex ante} so that in its place is some \(\{\hat{\pi}^\theta_n\}\sim\{p^\theta_n
    \}\) for economic pricing. Then, if the pair \(\{\hat{\pi}^\theta_n\}\sim\{p^\theta_n
    \}\) fails Proposition \ref{p1}, the \emph{ex post} (observed) risk-reward \(\{Y^{p
    }_n\}-\{Y^\pi_n\}\equiv\{Y^{p
    }_n\}-\{Y^{\hat{\pi}}_n\}+\{Y^{\hat{\pi}
    }_n\}-\{Y^\pi_n\}\) would be path-dependent with respect to \(\{\hat{\pi}^\theta_n\}\) and so \(\{{\pi}^\theta_n\}\), despite \emph{ex ante} economic risk-reward \(\{Y^{\hat{\pi}
    }_n\}-\{Y^\pi_n\}\) being path-independent by design. Thus, the risk-reward to be realised on average \emph{differs} after a round-trip in asset-prices even under identical agent-utility preference.
\end{remark}

For common tests, in which the null model (e.g. \(\theta_o=\{fair\text{-}coin\}\)) is objectively accurate, \(Q_{{\theta_o}}\equiv \mathbf{P}_{{\theta_o}}\), all above translate to requirements on the test models, \(Q_{\theta}\) and \(\hat{Q}_\theta\), in relation to the true \(\{B=\theta\}\)-law \(\mathbf{P}_{\theta}\). The condition of Proposition \ref{p1} (for path-independence) becomes a demand for equality, \(\hat{Q}_{\theta}=\mathbf{P}_\theta\) or \({Q}_\theta=\hat{Q}_\theta\). 

\section{Implications for Inferential Error}

Consider the error process \(\big|\{U[p^\theta_n]\}-\{U[\pi^\theta_n]\}\big|\) due to agent belief \(\{\pi^\theta_n\}\), where \(U\) is \emph{analytic} (e.g. for assessing a hazard-insurance), so that the object of interest is:\begin{equation}\label{Err}\{Err_n\}:=\big|\{p^\theta_n\}-\{\pi^\theta_n\}\big|\equiv\big|(\check{\pi}^\theta_n-\pi^\theta_n)+(p^\theta_n-\check{\pi}^\theta_n)\big|,\end{equation}where \(\check{\pi}^\theta_n\) is the would-be \emph{a posteriori} agent-belief given objectively accurate \emph{a priori} belief \(\check{\pi}^\theta_0\equiv p^\theta_0\), with the objective unconditional probability \(p^\theta_0\) generally unknown. Then, by definition: with the notation \(\sigma^{(\cdot)}_n:=({(\cdot)^\theta_n(\cdot)^{\theta_o}_n})^{\frac{1}{2}}\),\begin{align}
    \label{typ1}&
    \check{\pi}^\theta_n-\pi^\theta_n
    =(\rho^{\frac{1}{2}}-\rho^{-\frac{1}{2}})
    \sigma^{\check{\pi}}_n\sigma^\pi_n
    =\frac{(\rho^{\frac{1}{2}}-\rho^{-\frac{1}{2}})(\sigma^\pi_n)^2}{\rho^{\frac{1}{2}}\pi^\theta_n+\rho^{-\frac{1}{2}}\pi^{{\theta_o}}_n},\text{ \(\forall n\in\mathbb{N}\)},
\end{align}where \(\rho:=\frac{L[\check{\pi}^\theta_0]}{L[\pi^\theta_0]}\) is a \emph{conserved constant}: \(\frac{L[\check{\pi}^\theta_n]}{L[\pi^\theta_n]}\equiv\rho\) ((\ref{br1})).

Error (\ref{typ1}) is path-independent, and the only form that is, with respect to agent inference \(\{\pi^\theta_n\}\). It is \emph{a systematic bias} of the belief process \(\{\pi^\theta_n\}\), with 
a fixed degree \(\rho\) that \emph{does not} improve with data. For hazards \(\theta=+\) (e.g. disease), with \(p^+_0\ll\frac{1}{2}\), having \(\rho:=\frac{L[p^+_0]}{L[\pi^+_0]}>1\) 
means a bias for 
\(\theta_o=-\), {the status quo} of 'no hazards'. 

Error \(p^\theta_n-\check{\pi}^\theta_n\) ((\ref{Err})) due to \(Q_\theta\not=\mathbf{P}_\theta\) is path-dependent, and is independent of (\ref{typ1}). This error may or may not be systematic. 
For instance, in standard continuous settings its dynamic ((\ref{diff}), Item-\ref{ito'}, Appendix \ref{infps}) follows the difference of two diffusions with shared Wiener noise but different volatility, one representing the objective process and 
the other, the agent-modelled one. For i.i.d data ((\ref{diff}) having constant parameters), it is systematic and corresponds to persistent under- or overreactions to data.

The two error-components when combined can be mitigating or exacerbating  (e.g. status quo bias 'counters' overreaction), and their distinct and independent nature may be utilised in total-error management. 



\backmatter

\bmhead{Acknowledgement}Editor advice and reviewer comments are deeply appreciated.

\begin{appendices}
\section{Proof of Lemma \ref{l1}}\label{proofl1}
\begin{proof}
    Consider inferential odds \(\{L[{\pi}^{\theta}_{n}]\}\) and \(\{L[\hat{\pi}^{\theta}_n]\}\), \(L[{\pi}^{\theta}_{0}]=\alpha
    \) and \(L[\hat{\pi}^{\theta}_0]=\hat{\alpha}
    \), with 
    respective LR process \(\{L^{\theta}_n\}\) and \(\{\hat{L}^{\theta}_n\}\). Suppose there \(\exists g_n\) continuous, \(\forall n\in\mathbb{N}\), 
    such that \(L[\hat{\pi}^{\theta}_n]=g_n[L[{\pi}^{\theta}_{n}]]\). Then, by Bayes' Rule (\ref{br1}), for arbitrary \(n>1\) and at any \(m<n\in\mathbb{N}\),
    \begin{align}\label{sep}
    &\hat{\alpha}\hat{L}^{\theta}
        _{n}=g_n[\alpha L^{\theta}
        _{n}],\\
    \label{sep'}&\hat{L}^{\theta}
    _{n|m}=\frac{
g_n[\alpha L^{\theta}
_{m} L^{\theta}
_{n|m}]}{
g_m[\alpha L^{\theta}_m]}.\end{align}

Given any \(L^{\theta}_{m}=x_m\) non-null, 
 (\ref{sep'}) defines a mapping \(\hat{f}_{\{n,m,x_m;\alpha\}}\) from 
the \(\{x_m\}\)-conditional support \(\mathcal{L}_{n|m}^{x_m}\ni L^\theta_{n|m}
\) of one test to that 
of the other:\begin{equation}\label{sep''}
\hat{f}_{\{n,m,x_m;\alpha\}}:=\frac{g_n\circ \alpha x_m}{
g_m[\alpha x_m]}\text{ }\bigg|_{\mathcal{L}_{n|m}^{x_m}}
.\end{equation}

Recall (\ref{br1}) that LR processes become inferential odds iff. given \emph{a priori} odds, which are by definition independent of data and test models: mapping \(\{g_n\}\) must uphold (\ref{sep}-\ref{sep''}) 
\(\forall \alpha,\hat{\alpha}\in\mathbb{R}^+\). Immediately then, the domain of \(g_n\) is necessarily \(\mathbb{R}^+\), \(\forall n\in\mathbb{N}\).

Further, map (\ref{sep''}) must remain the same whether \(\alpha x_m=1\) or \(\alpha x_m=x\) arbitrary (keeping \(x_m\) fixed): thus \(\hat{f}_{\{n,m,x_m;\frac{1}{x_m}\}}=\hat{f}_{\{n,m,x_m;\frac{x}{x_m}\}}\), meaning,\begin{equation}\label{rawsep}\frac{g_n\circ x}{g_n}\bigg|_{\mathcal{L}_{n|m}^{x_m}}=\frac{g_m[x]}{g_m[1]}
,\text{ \(\forall x\in\mathbb{R}^+\)},\end{equation}where the RHS and so the LHS is invariant on \(\mathcal{L}_{n|m}^{x_m}\) for any given \(x\) (i.e. \(\alpha\)).

All above, (\ref{sep}-\ref{rawsep}), apply to any \(m'<n\) and \(x_{m'}\). Thus, provided \({\mathcal{L}_{n|m}^{x_m}\cap\mathcal{L}_{n|m'}^{x_{m'}}}\not=\emptyset\), which by Item-\ref{p.1} of Section \ref{technical} is indeed the case for any large enough \(n>m\&m'\), we have:\begin{equation}
        \label{m&m'}\frac{{g_m}[x]}{g_m[1]}=\frac{{g_{m'}}[x]}{g_{m'}[1]}
        ,\text{ \(\forall x\in\mathbb{R}^+\)}.
    \end{equation}Note that (\ref{sep}-\ref{m&m'}) apply to \(n\) arbitrarily large; so (\ref{m&m'}) holds \(\forall m\&m'\in\mathbb{N}\) finite. Then, set \(m'=1\) and write \(g_m=C_m\cdot g_1\), \(C_m\in\mathbb{R}^+\), \(m\in\mathbb{N}\), with \(C_1\equiv1\). Rewriting (\ref{rawsep}) yields:
    \begin{align}
\label{LL}&g_1\circ x\bigg|_{\mathcal{L}_{n|m}^{x_m}}=\frac{g_1[x]}{g_1[1]}\cdot g_1\bigg|_{\mathcal{L}_{n|m}^{x_m}},\text{ \(\forall x\in\mathbb{R}^+\)}.
\end{align}Equation (\ref{LL}) holds \(\forall n>m\in\mathbb{N}\), that is, it must hold over \({\mathcal{L}_{|m}^{x_m}}:=\bigcup_{n=m+1}^{\infty}
        \mathcal{L}^{x_m}_{n|m}\), the total \(\{L^{\theta}_{m}=x_m\}\)-conditional support, which by Item-\ref{p2} of Section \ref{technical} is dense in \(\mathbb{R}^+\).
        
        Then, under continuity (of the redundancy mapping by definition), equation (\ref{LL}) is in effect Cauchy: 
        \(g_1[xy]=\frac{{g_1[x]\cdot g_1[y]}}{g_1[1]},\text{ \(\forall x,y\in\mathbb{R}^+\)}\). 
Its only real continuous solution being \(g_1[\cdot]=(\cdot)^{\gamma}c\), \(\gamma\in\mathbb{R}\), \(c\in\mathbb{R}^+\), we have: by (\ref{sep}), \(\hat{L}^{\theta}_{1}=(\alpha L^{\theta}_{1})^{\gamma}c/{\hat{\alpha}}=(L^{\theta}_{1})^{\gamma}C\) given required independence from \({\alpha}\&{\hat\alpha}\) (i.e. \(c\equiv\frac{\hat{\alpha}}{\alpha^{\gamma}}C\)); the requirement of almost sure \emph{correct} resolution (Section \ref{st}) makes \(\gamma\) positive. 
\end{proof}

\begin{remark}\label{prfrmrk}
The proof makes no use of point-\ref{uni} of Section \ref{st}, hence Remark \ref{clarify}. Further, it concerns only arbitrary, distinct, moments; so its conclusion applies in continuous time. 
\end{remark}

\section{Common Sequential Inference}\label{infps}
Consider the log-LR process \(\{l^{\theta}_n\}:=\{\log L^{\theta}_n\}\) of a common test of \(B\)-outcomes \(\{\theta,{\theta_o}\}\) (as set up in Section \ref{st}). At any \(n\in\mathbb{N}\), the \(n\)th data-point \(D_n[n]\in V\) incurs:
\begin{align}
    &\Delta l^{\theta}_n(\cdot):=\log L^{\theta}_{n|n-1}(\cdot)=\log\frac{dQ_\theta|_{n}(\cdot|\mathcal{F}_{n-1})}{dQ_{{\theta_o}}|_{n}(\cdot|\mathcal{F}_{n-1})},\text{ with,}\\
    \label{frmdef}&{E}_{Q_\theta}[\Delta l^{\theta}_n|{\mathcal{F}_m}](B)=\frac{(-1)^{\mathbf{1}_{\{{\theta_o}\}}(B)}}{2}{E}_{Q_\theta}[(\Delta l^{\theta}_n)^2|{\mathcal{F}_m}](B),\text{ }\forall m<n\in\mathbb{N},
    \end{align}under 
    \emph{small increments}. 
    Thus i.i.d data 
guarantee \(B\)-detection as \(n\to\infty\) ((\ref{resolving})). 

\begin{enumerate}
\item{\label{takinglimit}\emph{Continuous-time inference.} In usual settings, for Wiener data, the above boil down to \emph{Wiener-drift testing} and become diffusions (
\citet{PShir}). For i.i.d data, time-homogeneity holds. 
Time-varying dataflow is modelled via 
\emph{time-change} 
so that the 
test dynamic acquires continuous time-varying parameters:
\begin{align}\label{w'}&dl^{\theta}_t(B)=(-1)^{\mathbf{1}_{\{{\theta_o}\}}(B)}\frac{(\sigma^{l}_t)^2}{2}
dt + \sigma^{l}_tdw_t,
\end{align}where \(\{dw_t\}\) is a standard Wiener process and \(\{\sigma^l_t\}\) the {signal-to-noise} of data (drift-difference-per-unit-volatility); 
note \(\lim_{t\to \infty}|l^{\theta}_t|=\infty\) for resolving tests.}

\item{\label{ito}\emph{Redundancy between Wiener-drift tests.}
Our results apply all the same (Remark \ref{prfrmrk}), but Ito's Lemma provides a sanity-check: any 2-differentiable time-homogeneous map \(
\hat{l}^{\theta}_t=g({l}^{\theta}_t)\) for tests obeying (\ref{w'}) and sharing the source of randomness \(\{dw_t\}\) must satisfy, under either \(B\)-outcome,\begin{equation}\label{ito1}g''= (-1)^{\mathbf{1}_{\{\theta_o\}}(B)}(g'-1)g';\end{equation}for divergent log-LRs it has only the identity as a solution, in either \(B\)-scenario.
 }

\item{\label{ito'}\emph{The path-dependent error component.}
In drift tests, differing test models lead to differing log-LRs (\ref{w'}). 
Given agent log-LR \(\{{l}^{\theta}_t\}\) with signal-to-noise \(\{\sigma^l_t\}\), %
and the objective \(\{\mathbf{l}^{\theta}_{t}\}\) with \(\{{\sigma}^\mathbf{l}_t\}\), 
their difference obeys, given \emph{shared} Wiener noise \(\{dw_t\}\):\begin{equation}\label{diff}\frac{d(\mathbf{l}^{\theta}_t-{l}^{\theta}_t)}{\sigma^\mathbf{l}_t+\sigma^{{l}}_t}(B)=(-1)^{\mathbf{1}_{\{{\theta_o}\}}(B)}\cdot\frac{\sigma^{\mathbf{l}}_t-\sigma^{{l}}_t}{2}
{dt} + \frac{\sigma^{\mathbf{l}}_t-\sigma^{{l}}_t}{\sigma^\mathbf{l}_t+\sigma^{{l}}_t}dw_t.\end{equation}}

\end{enumerate}



\end{appendices}



\begin{thebibliography}{00}


\bibitem[Breiman(1992)]{4.9}
Breiman, L. 1992. \textit{Probability}. SIAM, originally Addison-Wesley, Reading, Mass.

\bibitem[Bhatti \& Kim (2021)]{BK21}Bhatti, M.I. and Kim, J.H. 2021. {Towards a new paradigm for statistical evidence in the use of p-value}, \textit{Econometrics}, \textbf{9(1)}, 2, 1–3.

\bibitem[Bhatti \& Bodla (2008)]{BB08}Bhatti, M.I. and Bodla, M.A. 2008. {Empirical power comparison of non-nested tests for the EVM: Some Monte Carlo evidence}, \textit{International Journal of Applied Econometrics and Quantitative Studies}, \textbf{5(2)}.



\bibitem[{Cochrane(2005)}]{coh}
Cochrane, J.H. 2005. \textit{Asset Pricing}. Princeton University Press, Princeton.



\bibitem[Guidolin \& Timmermann(2007)]{GTim}
{Guidolin, M.} and {Timmermann, A.} 2007. {Properties of Equilibrium Asset Prices under Alternative Learning Schemes},
\textit{Journal of Economic Dynamics and Control} \textbf{31}
.




\bibitem[Jaiswal et al.
(2022)]{JCB}
{Jaiswal, S.}, {Chaturvedi, A.} and {Bhatti, M.I.} 2022. {Bayesian inference for unit root in smooth transition autoregressive models and its application to OECD countries},
\textit{Studies in Nonlinear Dynamics \& Econometrics} \textbf{26(1)}, 25-34.

\bibitem[Li et al.
(2026)]{lyl}
{Li, R.}, {Qin, Y.} and {Li, Y.} 2026. {Assessing estimation uncertainty
under model misspecification},
\textit{Scandinavian Journal of Statistics} \textbf{53(1)}, 268–290.


\bibitem[{Peskir and Shirayev(2006)}]{PShir}
{Peskir, G.} and {Shirayev, A.} 2006. \emph{
Optimal Stopping and Free-Boundary Problems. In: Lectures in Mathematics, ETH Zuerich}.
\newblock Birkhaeuser Verlag.





\end{thebibliography}

\end{document}